\documentclass[11pt]{amsart}
\usepackage{mathrsfs}
\usepackage{bbm}

\usepackage{amsmath}
\usepackage{amssymb}
\usepackage{amscd}
\usepackage{amsthm}

\newtheorem{thm}{Theorem}[section]
\newtheorem{lemma}[thm]{Lemma}
\newtheorem{prop}[thm]{Proposition}
\newtheorem{cor}[thm]{Corollary}

\newcommand{\bc}{\mathbb{C}}
\newcommand{\bz}{\mathbb{Z}}

\newcommand{\bp}{\mathbb{P}}
\newcommand{\py}{\mathbb{P}^1}
\newcommand{\pw}{\mathbb{P}^n}
\newcommand{\lpw}{L\mathbb{P}^n}

\newcommand{\fc}{\mathscr{C}}

\newcommand{\ch}{\mathscr{H}}
\newcommand{\chn}{\mathscr{H}(L\mathbb{P}^n)}
\newcommand{\chy}{\mathscr{H}(L\mathbb{P}^1)}

\newcommand{\htz}{H^2(L\mathbb{P}^n, \mathbb{Z})}
\newcommand{\htpoz}{H^2(\mathbb{P}^1, \mathbb{Z})}
\newcommand{\htpnz}{H^2(\mathbb{P}^n, \mathbb{Z})}

\newcommand{\hm}{f_{\phi}}
\newcommand{\hmlm}{f_{\phi, \,LM}}
\newcommand{\hmlg}{f_{\phi, \, LG}}
\newcommand{\hmlpw}{f_{\phi, \, L\mathbb{P}^n}}
\newcommand{\hmlcn}{f_{\phi, \, L \mathbb{C}^n}}

\newcommand{\df}{m_f}

\newcommand{\tr}{\mathscr{T}(S^1)}
\newcommand{\trk}{\mathscr{T}_k(S^1)}
\newcommand{\trkp}{\mathscr{T}_{k, p}(S^1)}

\newcommand{\dsy}{\mathscr{D}(S^1)}
\newcommand{\dsyk}{\mathscr{D}_k(S^1)}
\newcommand{\dsykp}{\mathscr{D}_{k, p}(S^1)}
\newcommand{\dsyyp}{\mathscr{D}_{1, p}(S^1)}
\newcommand{\dsyi}{\mathscr{D}_{\infty}(S^1)}

\newcommand{\auto}{\mathrm{Aut}(\lpw)}
\newcommand{\autoy}{\mathrm{Aut}(L\mathbb{P}^1)}

\begin{document}

\title{Holomorphic automorphisms of the loop space of $\mathbb{P}^n$}

\author{Ning Zhang}

\address{School of Mathematical Sciences\\
Peking University\\
Beijing, 100871\\ P. R. China}

\email{nzhang@math.pku.edu.cn}

\keywords{Loop space, Holomorphic self-map, Automorphism group,
Projective space}

\thanks{The author is grateful to L. Lempert and the referee for their
very helpful comments. This research was partially supported by the
National Natural Science Foundation of China grant 10871002.}

\begin{abstract}
The loop space $L\mathbb{P}^n$ of the complex projective space
$\mathbb{P}^n$ consisting of all $C^k$ or Sobolev $W^{k, \, p}$ maps
$S^1 \to \mathbb{P}^n$ is an infinite dimensional complex manifold.
We identify a class of holomorphic self-maps of $L\mathbb{P}^n$,
including all automorphisms.
\end{abstract}

\subjclass[2010]{32H02, 58B12, 46G20, 58D15, 58C10}

\maketitle

\section{Introduction \label{intro}}

Let $M$ be a finite dimensional complex manifold. We fix a
smoothness class $C^k$, $k=1,2, \cdots, \infty$, or Sobolev $W^{k,
\, p}$, $k=1,2, \cdots$, $1 \le p<\infty$, and consider the loop
space $LM=L_k M$, or $L_{k, \, p} M$, of all maps $S^1 \to M$ with
the given regularity. It carries a natural complex Banach/Fr\'echet
manifold structure (see \cite{l04}). In this paper, we identify a
class of holomorphic self-maps of the loop space $\lpw$ of the
complex projective space $\pw$. As a consequence, we compute the
group $\auto$ of holomorphic automorphisms of $\lpw$. This work was
directly motivated by \cite{mz, z03, lz04, z10}, in which certain
subgroups of $\autoy$ play a key role to study Dolbeault cohomology
groups with values in line bundles over $L\py$.

 There are two simple ways to construct holomorphic self-maps of a
given loop space $LM$. First, such a map can be obtained from a
family of holomorphic self-maps of $M$ smoothly parameterized by $t
\in S^1$. For example, let $G \simeq PGL(n+1, \bc)$ be the group of
holomorphic automorphisms of $\bp^n$. Its loop space $LG$ with
pointwise group operation is again a complex Lie group and acts on
$L\pw$ holomorphically; thus any element of $LG$ can be considered
as a holomorphic automorphism of $L\pw$. Second, let $\tr=\trk$
resp. $\trkp$ be the space of maps $\phi: S^1 \to S^1$ with the
following properties: for any $x$ in $L\bc=L_k\bc$ resp. $L_{k,
p}\bc$, the pull back $\phi^{\ast} x=x \circ \phi$ is still in
$L\bc=L_k\bc$ resp. $L_{k, p}\bc$, and the complex linear operator
$L\bc \ni x \mapsto x \circ \phi \in L\bc$ is continuous. If we set
the loop $x$ above to be the inclusion $S^1 \to \bc$, we see that
any $\phi$ in $\trk$ resp. $\trkp$ is a $C^k$ resp. $W^{k, p}$ map.
It is easy to verify that $\trk=C^k(S^1, S^1)$. Now any $\phi \in
\tr$ induces a holomorphic map
\begin{equation*}
    \hmlm: LM \ni x \mapsto x \circ \phi \in LM
\end{equation*}
(see Proposition \ref{hmlm}). We shall write $\hm$ for $\hmlpw$.

Recall that $\htz \simeq \bz$ (see \cite[Part II, Proposition
15.33]{cj} and \cite[Theorem 13.14]{p}). Any holomorphic self-map
$f$ of $L\pw$ induces a homomorphism $[\omega] \mapsto \df
[\omega]$, where $[\omega] \in \htz$, and $\df$ is a non-negative
integer. Our main result is

\begin{thm} \label{main}
Let $f$ be a holomorphic self-map of $L\pw$. Then $\df=1$ if and
only if there exist $\gamma \in LG$ and
  $\phi \in \tr$ such that
  \begin{equation} \label{decomposition}
     f=\gamma \circ \hm.
  \end{equation}
\end{thm}

Note that all constant loops are fixed points of $\hm$, and $\gamma$
does not share this property unless it is the identity. Hence the
decomposition of $f$ as in (\ref{decomposition}) is unique.

It is straightforward to verify that
\begin{equation} \label{relation}
\hm \circ \gamma =\hmlg(\gamma) \circ \hm, \hspace{4mm} \gamma \in
LG, \hspace{1mm} \phi \in \tr.
\end{equation}
Let $\dsy=\dsyk$ resp. $\dsykp$ be the space of bijections $\phi:
S^1 \to S^1$ such that both $\phi$ and $\phi^{-1}$ are in $\tr=\trk$
resp. $\trkp$. Then $\dsyk$ is the space of $C^k$ diffeomorphisms of
$S^1$. If $k=2, 3, \cdots$, or $k=p=1$, then $\dsykp$ is the space
of bijections $\phi: S^1 \to S^1$ such that $\phi, \phi^{-1} \in
W^{k, p}(S^1, S^1)$; and $\dsyyp$, where $p>1$, is the space of
bi-Lipschitz maps (see \cite[Theorem 4]{v} and \cite[Corollary
20.5]{hs}). Note that $\hm \in \auto$ if and only if $\phi \in
\dsy$, and $\dsy$ can be considered as a subgroup of $\auto$.

Combination of Theorem \ref{main} and (\ref{relation}) gives

\begin{cor} \label{semidirect}
The group $\auto$ is the semidirect product $LG \rtimes \dsy$.
\end{cor}

The group of holomorphic automorphisms of a compact complex manifold
is a complex Lie group. We would like to know whether $\auto$ can be
endowed with a natural complex Lie group structure. Recall that $LG$
is a complex Lie group. If $\dsy$ can be endowed with a manifold
structure such that $\auto=LG \rtimes \dsy$ becomes a complex Lie
group, then $\dsy \simeq \auto/LG$ is a complex Lie group. The space
$\dsyk$ can be considered as the open subset of $L_k S^1$ consisting
of embedded loops. With this manifold structure, $\dsyi$ is a Lie
group (and $\dsyk$, where $k<\infty$, is not a Lie group). It
follows from \cite[Proposition 3.3.2]{ps} that we cannot endow
$\dsy$ with a complex Lie group structure such that the inclusion
$\dsyi \to \dsy$ is a Lie group homomorphism.

Let $\fc=\fc(\lpw)$ be the set $\{\mu(\bp^1)\}$ of curves  in
$\lpw$, where $\mu$ ranges over all holomorphic embeddings $\bp^1
\to \lpw$ such that $\mu^{\ast}:$ $\htz$ $\to$ $\htpoz$ is an
isomorphism, and let $\ch=\chn$ be the space of holomorphic
self-maps $f$ of $\lpw$ with the following property: for any curve
$C \in \fc$, there exists a curve $C^{\,\prime} \in \fc$ such that
$f(C) \subset C^{\,\prime}$. It turns out that any holomorphic
self-map $f$ of $\lpw$ with $\df=1$ is in $\ch$ (see Section 3). To
prove Theorem \ref{main}, we only need to study maps in $\ch$.

This paper is organized as follows. In Section 2, we recall some
relevant facts about loop spaces and prove two propositions which
will be needed later on. Section 3 contains a complete proof of
Theorem \ref{main}. In the final Section 4, we try to classify maps
in $\ch$.

We refer to \cite{d, h} for the fundamentals of infinite dimensional
holomorphy.
\section{Preliminaries}

Let $\Psi: M \to M'$ be a holomorphic map between finite dimensional
complex manifolds. Then $L\Psi: LM \ni x \mapsto \Psi \circ x \in
LM'$ is holomorphic, and $L$ is functorial. For any $t \in S^1$, the
evaluation map $E_t: LM \ni x \mapsto x(t) \in M$ is holomorphic
(see \cite{l04}). The constant loops form a submanifold of $LM$,
which can be identified with $M$. Note that all elements of $LM$ can
be represented by absolutely continuous maps.

Let $G \simeq PGL(n+1, \bc)$ be as in Section 1. Applying the
functor $L$ to the holomorphic action $G \times \bp^n \to \bp^n$, we
obtain a holomorphic action $LG \times L\bp^n \to L\bp^n$.

\begin{prop} \label{hmlm}
Suppose $M$ is an $n$-dimensional complex manifold, where $0< n
<\infty$, and $\phi: S^1 \to S^1$ is a map. Then $\phi \in \tr$ if
and only if for any $x \in LM$, $x \circ \phi$ is still in $LM$, and
the map $\hmlm: LM \ni x \mapsto x \circ \phi \in LM$ is
holomorphic.
\end{prop}
\begin{proof}
First we show that if $\phi \in \tr$ and $x \in LM$, then $x \circ
\phi \in LM$. Let $(U, \Phi)$ be a coordinate chart of $M$, where $U
\subset M$ is open, and $\Phi$ is a biholomorphic map from $U$ to an
open subset of $\bc^n$. Then $LU \subset LM$ is open, and $L\Phi$ is
a biholomorphic map from $LU$ to an open subset of $L\bc^n$. If $x
\in LU$, then
\begin{equation*}
    x \circ \phi=\left(L\Phi\right)^{-1} \circ \hmlcn \circ
    L\Phi(x) \in LU \subset LM.
\end{equation*}
For a general loop $x \in LM$ and any $t_0 \in S^1$, there exist a
neighborhood $V' \subset S^1$ of $\phi(t_0)$, a coordinate chart
$(U, \Phi)$ of $M$ and $\tilde{x} \in LU$ such that
$x(t)=\tilde{x}(t)$, $t \in V'$. Note that $\phi$ is continuous.
Choose a neighborhood $V \subset S^1$ of $t_0$ such that $\phi(V)
\subset V'$, then $x \circ \phi(t)=\tilde{x} \circ \phi(t)$, $t \in
V$. So $x \circ \phi$ is $C^k$ resp. $W^{k, \, p}$, i.e. $x \circ
\phi \in LM$, and the map $\hmlm$ is well-defined.

Next we investigate the relationship between maps $\hmlm$ and
$\hmlcn$. Recall the complex structure of $LM$ as constructed in
\cite[Subsection 1.1]{ls}, where an open neighborhood of $y \in LM$
is mapped by the local chart $\varphi_y$ to an open subset of
$C^k(y^{\ast}TM)$ resp. $W^{k, \, p}(y^{\ast}TM)$, the space of
$C^k$ resp. $W^{k, \, p}$ sections of the pull back of the tangent
bundle of $M$ by $y$. It is straightforward to verify that
$\varphi_{y \circ \phi} \circ \hmlm \circ \varphi_y^{-1}$ is
precisely (the restriction of) the pull back
\begin{equation*}
    \phi^{\ast}: C^k(y^{\ast}TM) \to
    C^k(\phi^{\ast}y^{\ast}TM) \hspace{1mm}
    \mathrm{resp.} \hspace{0.8mm} W^{k, \, p}(y^{\ast}TM) \to
    W^{k, \, p}(\phi^{\ast}y^{\ast}TM).
\end{equation*}
The bundle $y^{\ast}TM$ over $S^1$ is always trivial, and the above
map can be considered as $\hmlcn$. Now we have $\varphi_{y \circ
\phi} \circ \hmlm \circ \varphi_y^{-1}=\hmlcn$. Thus $\phi \in \tr$
if and only if the map $\hmlm$ is well-defined and holomorphic.
\end{proof}

Let $i: \pw \to \lpw$ be the inclusion of the submanifold of
constant loops. Since $E_t \circ \hm \circ i$ is the identity, the
induced maps $i^{\ast}:$ $\htz$ $\to$ $\htpnz$, $\hm^{\ast}: \htz
\to \htz$ and $E_t^{\ast}:$ $\htpnz$ $\to$ $H^ 2(\lpw, \bz)$ are all
isomorphisms.

The following simple proposition will be very useful:

\begin{prop} \label{embedding}
Let $\mu: \py \to \lpw$ and $\tau: \pw \to \lpw$ be holomorphic maps
such that $\mu^{\ast}: \htz \to H^2(\py, \bz)$ and $\tau^{\ast}:$
$\htz$ $\to$  $\htpnz$ are isomorphisms. Then:
\begin{enumerate}
  \item[\textup{(a)}] For any $t \in S^1$, the map $E_t
  \circ \mu: \py \to \pw$ is an embedding whose image is a
  projective line. In particular, $\mu$ itself is an embedding.

  \item[\textup{(b)}] There exists
  $\gamma \in LG$ (considered as an automorphism of $\lpw$)
  such that $\tau=\gamma|_{\pw}$.
\end{enumerate}
\end{prop}
\begin{proof}
(a) Note that $(E_t \circ \mu)^{\ast}=\mu^{\ast} \circ E_t^{\ast}:
\htpnz \to H^2(\py, \bz)$ is an isomorphism, which maps the first
Chern class of the hyperplane section bundle over $\pw$ to the first
Chern class of the hyperplane section bundle over $\py$. Therefore
the inverse image of any hyperplane of $\pw$ under $E_t \circ \mu$
is either a hyperplane (i.e. a single point) or the entire $\py$.
The proof will be by induction on $n$. If $n=1$, then $E_t \circ
\mu$ is conformal. If $n> 1$, take a hyperplane $H \subset \pw$
containing two different points in $E_t \circ \mu(\py)$. Then we
must have $E_t \circ \mu(\py) \subset H \simeq \mathbb{P}^{n-1}$.

(b) The restriction of $\tau$ to any projective line in $\pw$ can be
considered as the map $\mu$ . In view of part (a), for any $t \in
S^1$, $E_t \circ \tau: \pw \to \pw$ is injective; hence $E_t \circ
\tau \in G$. Define $\gamma: S^1 \ni t \mapsto E_t \circ \tau \in
G$. Then
\begin{equation} \label{auto}
    \gamma(t)(\zeta)=\tau(\zeta)(t), \hspace{3mm} \zeta \in \pw, t
    \in S^1.
\end{equation}
Note that $\gamma$ is $C^k$ resp. $W^{k, \,p}$ (i.e. $\gamma \in
LG$) if $\gamma(t)(\zeta_j)$, $j=1, 2, \cdots, n+2$, are $C^k$ resp.
$W^{k, \, p}$ maps of $t$, where $\{\zeta_1, \zeta_2, \cdots,
\zeta_{n+2}\} \subset \pw$ is any given set of $n+2$ points in
general position. It follows from (\ref{auto}) that $\gamma \in LG$
and $\tau=\gamma|_{\pw}$.
\end{proof}

\section{Proof of Theorem \ref{main}}

We begin with two results concerning holomorphic self-maps of loop
spaces of the type $\hmlm$.

\begin{prop} \label{rotation}
Let $f$ be a holomorphic self-map of $LM$. Then $f=\hmlm$ for some
$\phi \in \tr$ if and only if for any $t \in S^1$, there exists $t'
\in S^1$ such that $E_t \circ f =E_{t'}$.
\end{prop}
\begin{proof}
If $f=\hmlm$, then $E_t \circ f=E_{\phi(t)}$. For the other
direction, define the map $\phi: S^1 \ni t \mapsto t' \in S^1$. Then
$f(x)=x \circ \phi$. It follows from Proposition \ref{hmlm} that
$\phi \in \tr$.
\end{proof}

\begin{lemma} \label{lcrotation}
Let $f$ be a self-map of $L\bc^n$. Then $f=\hmlcn$ for some $\phi
\in \tr$ if and only if $f$ is continuous complex linear and
\begin{equation} \label{subsets}
f(x)(S^1) \subset x(S^1)
\end{equation}
for all $x \in L\bc^n$.
\end{lemma}
\begin{proof}
One direction being trivial, we shall only verify the sufficiency
part of the claim. Let $e_1, e_2, \cdots, e_n$ be the standard basis
of $\bc^n \subset L\bc^n$. For any $x \in L\bc^n$, we write
$x=\sum_{j=1}^n x_je_j$, where $x_j \in L\bc$. In view of
(\ref{subsets}), we have $f(x_j e_j)=y_j e_j$, where $y_j \in L\bc$
and $y_j(S^1) \subset x_j(S^1)$. Thus $f$ induces continuous complex
linear maps $f_j: L\bc \ni x_j \mapsto y_j \in L\bc$, $j=1, 2,
\cdots, n$. Now $f(x)=\sum_{j=1}^n f_j(x_j) e_j$. Let
$\bc^{\ast}=\bc \setminus \{ 0\}$. Then $f_j(L\bc^{\ast}) \subset
L\bc^{\ast}$ and $f_j(1)=1$. It follows from \cite[Lemma 3.1]{z03}
that for any $t \in S^1$, there exists $t_j \in S^1$ such that $E_t
\circ f_j=E_{t_j}$. Therefore
\begin{equation} \label{eth}
E_t \circ f(x)=\sum_{j=1}^n x_j(t_j)e_j.
\end{equation}
Let $x_0 \in L\bc$ be an embedded loop. Setting $x=\sum_{j=1}^n x_0
e_j$ in (\ref{subsets}) and (\ref{eth}), we obtain that
$x_0(t_1)=x_0(t_2)=\cdots=x_0(t_n)$. Hence $t_1=t_2=\cdots=t_n$ and
$E_t \circ f=E_{t_1}$. By Proposition \ref{rotation}, $f=\hmlcn$ for
some $\phi \in \tr$.
\end{proof}

\begin{prop} \label{2points}
Let $x$ and $y$ be two different points of $\lpw$. Then there exists
a curve $C \in \fc$ through both $x$ and $y$ if and only if $x(t)
\not= y(t)$ for all $t \in S^1$.
\end{prop}
\begin{proof} The ``only if'' direction follows from Proposition
\ref{embedding}(a). To show the ``if'' direction, consider the
natural projection $\pi: \bc^{n+1} \setminus \{0\} \to \bp^n$ and
the holomorphic map $L\pi: L(\bc^{n+1} \setminus \{0\}) \to \lpw$.
Let $\tilde{x}, \tilde{y} \in L(\bc^{n+1} \setminus \{0\})$ be such
that $\tilde{x}(t)$ and $\tilde{y}(t)$ are linearly independent for
all $t \in S^1$. Then the map
\begin{equation} \label{construct}
\mu=\mu_{\tilde{x}, \,\tilde{y}}: \bp^1 \ni [Z_0, Z_1] \mapsto L\pi
\left(Z_0 \tilde{x}+Z_1 \tilde{y} \right) \in \lpw,
\end{equation}
where $Z_0, Z_1$ are homogeneous coordinates on $\bp^1$, is
well-defined and holomorphic. For any $t \in S^1$, the map $E_t
\circ \mu$ is an embedding whose image is a projective line in
$\bp^n$; thus $\mu$ is an embedding and $\mu^{\ast}: \htz \to
\htpoz$ is an isomorphism.

The map $L\pi$ is surjective (for $C^k$ loops see \cite[Lemma
2.2]{ls}, and the proof of \cite[Lemma 2.2]{ls} also works for
$W^{k, \,p}$ loops). Take $\tilde{x} \in (L\pi)^{-1}(x)$ and
$\tilde{y} \in (L\pi)^{-1}(y)$ in (\ref{construct}), then the image
of $\mu$ is a curve $C \in \fc$ through both $x$ and $y$.
\end{proof}

From now on, we shall concentrate on holomorphic self-maps $f$ of
$\lpw$, and curves in $\fc$ will be needed throughout the rest of
the paper.

\begin{prop} \label{degree0}
If $\df=0$, then $f$ is constant.
\end{prop}
\begin{proof}
If $\eta: \bp^1 \to \pw$ is a holomorphic map such that
$\eta^{\ast}: \htpnz \to \htpoz$ is zero, then $\eta$ is constant.
So for any curve $C \in \fc$ and any $t \in S^1$, the map $E_t \circ
f|_C: C \to \pw$ is constant; hence $f|_C$ is constant. It follows
from Proposition \ref{2points} that for any $x \in \lpw$ and any
constant loop $\zeta \in \pw \setminus x(S^1)$, there exists a curve
$C \in \fc$ through both $x$ and $\zeta$; thus $f$ is constant.
\end{proof}

Recall the mapping space $\ch$ as defined in Section \ref{intro}. If
$f \in \ch$, then the topological degree of $f|_C: C \to
C^{\,\prime}$ is $\df$. In particular, $f(C)=C^{\,\prime}$ if $\df
\ge 1$, and $f|_C$ is one-to-one if $\df=1$. By Proposition
\ref{embedding}(a), any holomorphic self-map $f$ of $\lpw$ with
$\df=1$ is in $\ch$. Hence $LG \subset \ch$, and $\hm \in \ch$ for
any $\phi \in \tr$.

\begin{prop} \label{image}
Suppose $f \in \ch$, $\df \ge 1$, $x \in \lpw$, and $\zeta \in \pw
\setminus x(S^1)$ is a constant loop. If either $\df=1$ or $f(x)
\not\in f(\pw)$, then $f(x)(t) \not= f(\zeta)(t)$ for all $t \in
S^1$.
\end{prop}
\begin{proof}
By Proposition \ref{2points}, there exists a curve $C \in \fc$
through both $x$ and $\zeta$. If either $\df=1$ or $f(x) \not\in
f(\pw)$, then $f(x) \not=f(\zeta)$. Since both $f(x)$ and $f(\zeta)$
are in $f(C) \in \fc$, it follows from Proposition \ref{2points}
that $f(x)(t) \not= f(\zeta)(t)$ for all $t \in S^1$.
\end{proof}

\begin{thm} \label{lpnrotation}
Let $f$ be a holomorphic self-map of $\lpw$. Then $f=\hm$ for some
$\phi \in \tr$ if and only if every constant loop is a fixed point
of $f$.
\end{thm}
\begin{proof}
The necessity is obvious. Regarding sufficiency, note that
$f|_{\pw}$ is the identity; hence $\df=1$ and $f \in \ch$. It
follows from Proposition \ref{image} that
\begin{equation} \label{subsets2}
f(x)(S^1) \subset x(S^1)
\end{equation}
for all $x \in \lpw$. Let $$U_0=\{ [Z_0, Z_1, \cdots, Z_n] \in \pw:
Z_0 \not=0\},$$ where $Z_0, Z_1, \cdots, Z_n$ are homogeneous
coordinates on $\pw$. Now consider $U_0$ as $\bc^n$. It follows from
(\ref{subsets2}) that $f(L\bc^n) \subset L\bc^n$. Next we show that
$f|_{L\bc^n}$ is complex linear.

Let $x=(x_1, \cdots, x_n) \in L\bc^n$ and $y=(y_1, \cdots, y_n) \in
L(\bc^n \setminus \{0\})$. If we choose $\tilde{x}=(1, x_1, \cdots,
x_n)$ and $\tilde{y}=(0, y_1, \cdots, y_n)$ in (\ref{construct}),
then the image of $\mu$ is a curve $C_{x, \, y} \in \fc$, and
\begin{equation*}
    C_{x, \, y} \cap
    L\bc^n=C_{x, \, y} \setminus \{L\pi(\tilde{y}) \}=\{x+\lambda y \in L\bc^n:
    \lambda \in \bc\},
\end{equation*}
where $L\pi(\tilde{y}) \in L(\pw \setminus \bc^n)$. Note that
$f(C_{x, \, y})$ is also a curve in $\fc$, and $f$ maps $C_{x, \,
y}$ conformally onto $f(C_{x, \, y})$. By Proposition
\ref{embedding}(a), for any $t \in S^1$, $E_t \circ f(C_{x, \, y})$
is a projective line in $\pw$. In view of (\ref{subsets2}), we have
\begin{equation} \label{affineline}
    E_t \circ f(C_{x, \, y}) \cap \bc^n=\{E_t \circ f(x+\lambda y):
    \lambda \in \bc\}.
\end{equation}
The above set must be an affine line of $\bc^n$. As a function of
$\lambda$, $E_t \circ f(x+\lambda y)$ maps $\bc$ bijectively onto
the affine line in (\ref{affineline}); therefore it is a polynomial
of degree one. So
\begin{equation*}
    E_t \circ f(x+\lambda y)=\left[f(x+y)(t)-f(x)(t)
    \right]\lambda+f(x)(t)
\end{equation*}
for all $t \in S^1$. Thus $$f(x+\lambda y)=
\left[f(x+y)-f(x)\right]\lambda+f(x),$$ i.e. $\bc \ni \lambda
\mapsto f(x+\lambda y) \in L\bc^n$ is a polynomial of degree one,
where $x \in L\bc^n$ and $y \in L(\bc^n \setminus \{0\})$. Since
$L(\bc^n \setminus \{0\})$ is dense in $L\bc^n$, $f(x+\lambda y)$ is
a polynomial of $\lambda$ of degree less than or equal to one for
all $x, y \in L\bc^n$. Hence $f|_{L\bc^n}$ is a polynomial of degree
one (see \cite[Section 2.2]{h}). As $f(0)=0$, $f|_{L\bc^n}$ is
complex linear.

By (\ref{subsets2}) and Lemma \ref{lcrotation}, we have
$f|_{L\bc^n}=\hmlcn$ for some $\phi \in \tr$; thus $f=\hm$ on the
connected manifold $\lpw$.
\end{proof}

\noindent {\it Proof of Theorem \ref{main}.} The sufficiency is
obvious. Regarding necessity, by Proposition \ref{embedding}(b),
there exists $\gamma \in LG$ such that $f|_{\bp^n}=\gamma|_{\bp^n}$.
Then every constant loop is a fixed point of $\gamma^{-1} \circ f$.
It follows from Theorem \ref{lpnrotation} that $\gamma^{-1} \circ
f=\hm$ for some $\phi \in \tr$, i.e. $f=\gamma \circ \hm$. \qed
\section{The mapping space $\ch$}

In this section, we continue to study maps in $\ch$. Recall that all
holomorphic self-maps $f$ of $\lpw$ with $\df \le 1$ are in $\ch$.

\begin{thm} \label{converse}
Let $f \in \ch=\chn$.
\begin{enumerate}
\item[(a)] If $n \ge 2$, then $\df \le 1$.

\item[(b)] If $n=1$ and $\df \ge 2$, then $f(L\py)=f(\py)$
(i.e. $f(C)=f(\py)$ for all $C \in \fc$).
\end{enumerate}
\end{thm}
\begin{proof}
(a) If $f$ is constant, then $\df=0$. If $f$ is non-constant, by
Proposition \ref{degree0} we have $\df \ge 1$. Fix $t \in S^1$ and
define $\alpha=E_t \circ f|_{\pw}: \pw \to \pw$. Then $\alpha$
induces a homomorphism $[\omega] \mapsto \df [\omega]$, where
$[\omega] \in \htpnz$. So the topological degree of $\alpha$ is
$\df^n \not=0$. Let $\zeta^{\ast} \in \pw$ be a regular value of
$\alpha$. Take $\zeta \in \alpha^{-1}(\zeta^{\ast})$ and choose a
hyperplane $H \subset \pw$ such that $\zeta \not\in H$. For any $w
\in H$, the projective line $P_{\zeta, \,w}$ through both $\zeta$
and $w$ is in $\fc$; thus $f(P_{\zeta, \,w}) \in \fc$. In view of
Proposition \ref{embedding}(a), $\alpha(P_{\zeta, \,w})$ is a
projective line of $\pw$. The topological degree of the map
$P_{\zeta, \,w} \stackrel{\alpha}{\to} \alpha(P_{\zeta, \,w})$ is
$\df$. If $\df \ge 2$, then $P_{\zeta, \,w} \setminus \{\zeta\}$
must contain at least one point in $\alpha^{-1}(\zeta^{\ast})$. For
different $w \in H$, the sets $P_{\zeta, \,w} \setminus \{\zeta\}$
are disjoint; hence $\alpha^{-1}(\zeta^{\ast})$ is not a finite set,
which is a contradiction.

(b) Note that $f(\bp^1) \in \fc$. By Proposition \ref{embedding}(b)
(in which we choose $\tau$ to be a suitable holomorphic embedding
$\bp^1 \to L\bp^1$ with image $f(\bp^1)$), there exists $\gamma \in
LG$ such that $\gamma(\bp^1)=f(\bp^1)$; thus $\gamma^{-1} \circ
f(\bp^1)= \bp^1$. Without loss of generality, we may assume that
$f(\bp^1)= \bp^1$.

Let $\rho=f|_{\bp^1}$. Then $\rho: \bp^1 \to \bp^1$ is an
$\df$-sheeted branched covering. Choose a non-empty open subset $W
\subset \bp^1$ such that $\rho|_W$ is one-to-one. Then $LW \subset
L\bp^1$ is a non-empty open subset. We claim that
\begin{equation} \label{constant}
f(x) \in \bp^1 \hspace{1mm} \mathrm{for} \hspace{1.2mm} x \in LW.
\end{equation}
Otherwise there would exist $x_0 \in LW$ such that $y_0=f(x_0)
\not\in \bp^1$. The set $y_0(S^1)$ is not finite; thus we can find a
regular value $\zeta^{\ast}$ of $\rho$ in $y_0(S^1)$. Take $\zeta_1,
\zeta_2 \in \rho^{-1}(\zeta^{\ast})$, $\zeta_1 \not= \zeta_2$. It
follows from Proposition \ref{image} that $\zeta_1, \zeta_2 \in
x_0(S^1) \subset W$; then $\rho|_W$ is not one-to-one, which is a
contradiction.

By (\ref{constant}), maps $E_t \circ f|_{LW}$ are independent of $t
\in S^1$; hence maps $E_t \circ f$ are independent of $t$ on the
connected manifold $L\bp^1$, i.e. $f(L\bp^1) \subset \bp^1$.
\end{proof}

The maps as in Theorem \ref{converse}(b) do exist: Let $\rho: \bp^1
\to \bp^1 \subset L\bp^1$ be a holomorphic map with topological
degree $m \ge 2$, $t \in S^1$ and $\gamma \in LG$. Then $f=\gamma
\circ \rho \circ E_t \in \chy$ and $\df=m$.

\end{document}